\documentclass[12pt,letterpaper]{amsart}
\usepackage{amsmath,amssymb,amsfonts,amsthm}
\usepackage{mathrsfs}
\usepackage[all]{xy}
\usepackage{stmaryrd}
\usepackage{mathtools}
\usepackage{indentfirst}
\usepackage{verbatim}
\usepackage{hyperref}
\usepackage{color}
\usepackage[usenames,dvipsnames]{xcolor}
\usepackage{ulem}
\usepackage[shortlabels]{enumitem}

\title{Fujita's Conjecture for Quasi-Elliptic Surfaces}
\author{Yen-An Chen}
\subjclass[2010]{Primary 14C20, Secondary 14G17}
\thanks{The author was partially supported by NSF research grants no: DMS-1801851, DMS-1840190 and by a grant from the Simons Foundation; Award Number: 256202. }
\address{Department of Mathematics, University of Utah, Salt Lake City, UT 84112, USA}
\email{yachen@math.utah.edu}

\newtheorem{thm}{Theorem}[section]
\newtheorem{prop}[thm]{Proposition}
\newtheorem{cor}[thm]{Corollary}
\newtheorem{lem}[thm]{Lemma}

\newtheorem{conj}[thm]{Conjecture}
\theoremstyle{definition}
\newtheorem{defn}[thm]{Definition}
\newtheorem*{acks}{Acknowledgements}
\theoremstyle{remark}
\newtheorem{rmk}[thm]{Remark}
\newtheorem*{pf}{Proof}

\newcommand\cE{{\mathcal{E}}}

\newcommand\cI{{\mathcal{I}}}

\newcommand\cO{{\mathcal{O}}}

\newcommand\bP{{\mathbb P}}
\newcommand\bQ{{\mathbb Q}}

\newcommand\rw{\rightarrow}
\newcommand\wt[1]{\widetilde{#1}}

\begin{document}

\maketitle

\begin{abstract}
We show that Fujita’s conjecture is true for quasi-elliptic surfaces. Explicitly, for any quasi-elliptic surface $X$ and an ample line bundle $A$ on $X$, we have $K_X+tA$ is base point free for $t\geq 3$ and is very ample for $t\geq 4$.
\end{abstract}

\section{Introduction}\label{sec1}
Let $X$ be a smooth projective variety and $A$ be an ample line bundle. 
One of the crucial classical problems is understanding under what conditions the adjoint linear system $K_X + A$ is base point free or very ample. 
Thanks to Serre's theorem, we know that $K_X + tA$ is very ample for $t$ sufficiently large, and there is great interest in understanding the smallest value of $t$ for which this holds.
The following conjecture is due to Fujita in \cite{fujita1988problems}. 

\begin{conj}[Fujita]
Let $X$ be a smooth projective variety of dimension $n$ and $A$ be an ample line bundle. 
Then $K_X+tA$ is base point free (resp. very ample) whenever $t\geq n+1$ (resp. $t\geq n+2$).
\end{conj}

This conjecture for curves follows from the Riemann-Roch theorem. 
In characteristic zero, the conjecture is entirely proved for surfaces by Reider's theorem \cite{reider1988vector}. 
For the base point freeness part of this conjecture in characteristic zero, it has been proved up to dimension five in \cite{ein1993global}, \cite{helmke1997fujita}, \cite{kawamata1997fujita}, and \cite{ye2020fujita}. 
In positive characteristic, Shepherd-Barron showed in \cite{shepherd1991unstable} that the conjecture is true for surfaces that are neither quasi-elliptic (see Definition~\ref{defn_quasi_ell}) nor of general type. Recently, Gu, Zhang, and Zhang claimed in \cite{gu2020counterexamples} that there are counterexamples to the conjecture for surfaces. Those examples are of general type. 

In this paper, we show that Fujita's conjecture is true for quasi-elliptic surfaces. 
\begin{thm}[=Theorem~\ref{main}]
Fujita's conjecture is true for quasi-elliptic surfaces $X$. 
(See Definition~\ref{defn_quasi_ell}.)
That is, given a quasi-elliptic surface $X$ and any ample line bundle $A$ on $X$, we have 
\begin{enumerate}
\item $K_X+tA$ is base point free for $t\geq 3$; and 
\item $K_X+tA$ is very ample for $t\geq 4$. 
\end{enumerate}
\end{thm}

To prove this result, we follow the ideas of \cite{di2015effective} and make a careful case-by-case study.  
Note that, in \cite{di2015effective}, it is proved that, when $p=3$, $K_X+tA$ is base point free for $t\geq 4$ and it is very ample for $t\geq 8$; 
and when $p=2$, $K_X+tA$ is base point free for $t\geq 5$ and it is very ample for $t\geq 19$. 

\begin{acks}
The author would like to thank his advisor Christopher D. Hacon for many helpful suggestions and encouragement. 
\end{acks}

\section{Preliminaries}\label{sec2}
In this section, we recall some definitions and results which will be used later. 
We will always assume the base field $k$ is algebraically closed and of positive characteristic $p$. 

\begin{lem}\label{Hodeg_ineq}
Let $X$ be a smooth projective surface over $k$ and $N$ a nef divisor on $X$. 
Then for any divisor $D$ on $X$, we have 
\[N^2D^2\leq (N.D)^2.\] 
Moreover, if $N$ is ample, then the equality holds only when $D$ is numerically proportional to $N$. 
\end{lem}
\begin{pf}
Since we can approximate nef divisors by ample $\bQ$-divisors and the desired inequality is homogeneous, 
we can reduce to the case when $N$ is ample.

Now we consider $E = (N.D)N - N^2D$. 
Notice that $E.N = 0$. 
Then, by the Hodge index theorem, we have $E^2\leq 0$, and we get the desired inequality. 
Moreover, the equality holds only when $E\equiv 0$, that is, $D$ is numerically proportional to $N$. \qed
\end{pf}

\begin{defn}\label{defn_quasi_ell}
A smooth projective surface $X$ over $k$ is said to be {\it quasi-elliptic} if there is a fibration $f : X \rw C$ where $C$ is a smooth curve such that $f_*\cO_X = \cO_C$ and such that the general fibers of $f$ 
are rational curves with one (ordinary) cusp. 
Such $f:X\rw C$ is called a quasi-elliptic fibration. 
\end{defn}
\begin{rmk}
The general fibers of $f$ have arithmetic genus 1. 
Moreover, by a result of Tate in \cite{tate1952genus}, quasi-elliptic surfaces exist only when $p=2$ or $3$. 
\end{rmk}

\begin{defn}[\cite{shepherd1991unstable},\cite{di2015effective}]
A rank-two vector bundle $\cE$ on $X$ is {\it unstable} if it fits into a short exact sequence, which will be called a de-stabilizing sequence for $\cE$,  
\[\xymatrix{0\ar[r] & \cO_X(D_1)\ar[r] & \cE \ar[r] & \cI_Z\cdot\cO_X(D_2) \ar[r] & 0}\]
where $D_1$ and $D_2$ are Cartier divisors, $I_Z$ is the ideal sheaf of a finite subscheme $Z$ of $X$, and $D_1 - D_2\in C_{++}(X)$, the positive cone of $\textnormal{NS}(X)$. 
Notice that $Z$ could be empty, and by convention, $\cI_Z = \cO_X$ when $Z$ is empty. 
We also recall that 
\[C_{++}(X) = \{x\in\textnormal{NS}(X)\vert x^2>0 \mbox{ and } x.H>0 \mbox{ for some ample divisor } H\} = \{x\in\textnormal{NS}(X)\vert x^2>0 \mbox{ and } x \mbox{ is big}\}.\] 
\end{defn}

\begin{defn}[\cite{di2015effective}]\label{m_unstable}
A big divisor $D$ on a smooth surface $X$ with $D^2>0$ is {\it $m$-unstable} for a positive integer $m$ 
if $h^1(X,\cO_X(-D)) \neq 0$ and there exists a nonzero effective divisor $E$ such that 
$mD-2E$ is big and $(mD-E).E\leq 0$.
\end{defn}

In \cite{bogomolov1978holomorphic}, Bogomolov showed that, in characteristic zero, every rank-two vector bundle $\cE$ on a smooth surface 
with $c_1^2(\cE)>4c_2(\cE)$ is unstable. 
Also, in positive characteristic, there is a result related to the unstability of vector bundles. 

\begin{thm}[{\cite[Theorem 1]{shepherd1991unstable}}]\label{lift}
Let $\cE$ be a rank-two vector bundle on a smooth projective surface $X$ over an algebraically closed field $k$ of positive characteristic $p>0$ 
such that $c_1^2(\cE)>4c_2(\cE)$. 
Then there exists an integral surface $Y$ contained in the ruled threefold $\bP(\cE)$ such that 
\begin{enumerate}
\item the composition $\rho : Y \to X$ is purely inseparable of degree $p^e$ for some $e>0$; and
\item $(F^e)^*\cE$ is unstable where $F : X \rw X$ is the absolute Frobenius morphism. 
\end{enumerate}
Moreover, we have \[K_Y \equiv \rho^*\left(K_X-\frac{p^e-1}{p^e}(D_1-D_2)\right)\]
where \[\xymatrix{0\ar[r] & \cO_X(D_1)\ar[r] & F^{e*}\cE \ar[r] & \cI_Z\cdot\cO_X(D_2) \ar[r] & 0}\]
is a de-stablizing sequence for $(F^e)^*\cE$. 
\end{thm}

We recall the construction of $Y$ in Theorem~\ref{lift}. 
Assume we have shown that $(F^n)^*\cE$ is unstable for some positive integer $n$. 
Let $e$ be the smallest one such that $\wt{\cE} := (F^e)^*\cE$ is unstable. 
We have the following cartesian diagram 
\[\xymatrix{\bP(\wt{\cE}) \ar[r]^-\psi \ar[d]_-{\wt{\pi}} & \bP(\cE) \ar[d]^-\pi \\
X \ar[r]_-{F^n} & X.}\]
From a de-stablizing sequence for $(F^e)^*\cE$, we have a surjection $(F^e)^*\cE \rw \cI_Z\cdot\cO_X(D_2)$, which gives a quasi-section $Y'\subset \bP(\wt{\cE})$. 
Then $Y$ is the schematic image of $Y'$ in $\bP(\cE)$. 

\begin{lem}\label{key_lem}
If $D$ is big with $D^2>0$ and $h^1(X,\cO_X(-D))\neq 0$, 
then $D$ is $p^e$-unstable for some $e>0$. 
\end{lem}
\begin{pf}
Indeed, this is contained in \cite[Lemma 16]{shepherd1991unstable}. 
The reader can also see \cite[Remark 2.10]{di2015effective}. 
For the reader's convenience,  we include the proof. 

Since $h^1(X,\cO_X(-D))\neq 0$, there exists a non-split short exact sequence 
\[\xymatrix{0\ar[r] & \cO_X \ar[r] & \cE \ar[r] & \cO_X(D) \ar[r] & 0}\]
given by a nonzero element of $\mbox{Ext}^1(\cO_X(D),\cO_X) \cong H^1(X,\cO_X(-D))$, 
where $\cE$ is a vector bundle of rank two. 
Note that $c_1^2(\cE)-4c_2(\cE) = D^2>0$. 
By Theorem~\ref{lift}, we have the following diagram. 
\[\xymatrix{ & & 0 \ar[d] & & \\
 & & \cO_X \ar[d]^-{g_1} & & \\
0 \ar[r] & \cO_X(D_1) \ar[r]^-{f_1} \ar[rd]_-\tau & (F^e)^*\cE \ar[r]^-{f_2} \ar[d]^-{g_2} & \cI_Z\cdot\cO_X(D_2) \ar[r] & 0 \\
 & & \cO_X(p^eD) \ar[d] & & \\
 & & 0 & & }\]
We claim that $\tau = g_2\circ f_1$ is not zero. 
Indeed, if $\tau = 0$, then $f_1 = g_1\circ\tau'$ where $\tau'$ is a nonzero map from $\cO_X(D_1)$ to $\cO_X$. 
That means $-D_1$ is linearly equivalent to an effective divisor. 
Now notice that $D_1+D_2\equiv c_1((F^e)^*\cE) \equiv p^eD$ is big and, for any ample divisor $H$, we have 
\[0<p^eD.H = (D_1+D_2).H = -(D_1-D_2).H - (-2D_1).H < 0, \mbox{ which is impossible.}\] 

Hence, we have $\tau\neq 0$ and so $D_2 \equiv c_1((F^e)^*\cE) -D_1\equiv p^eD-D_1$ is effective. 
So $p^eD-2D_2\equiv D_1-D_2$ is big and 
\[(p^eD-D_2).D_2 = D_1.D_2 = c_2((F^e)^*\cE) - \deg Z = -\deg Z \leq 0.\] 

Also $D_2\neq 0$ since otherwise the vertical exact sequence 
\[\xymatrix{0\ar[r] & \cO_X \ar[r] & \cE \ar[r] & \cO_X(D) \ar[r] & 0}\]
splits, which is a contradiction.

To sum up, $D$ is $p^e$-unstable. \qed
\end{pf}

\begin{prop}\label{key_prop}
Let $\pi: Y \to X$ be a birational morphism between two smooth surfaces and let $\wt{D}$ be a big Cartier divisor on $Y$ such that $\wt{D}^2>0$. 
Assume there is a non-zero effective divisor $\wt{E}$ such that 
$\wt{D}-2\wt{E}$ is big and 
$(\wt{D}-\wt{E}).\wt{E}\leq 0$. 

Set $D = \pi_*\wt{D}$, $E = \pi_*\wt{E}$ and $\alpha = D^2-\wt{D}^2$. 
If $D$ is nef and $E$ is a nonzero effective divisor, then
$0\leq D.E<\alpha/2$ and 
$D.E-\alpha/4\leq E^2\leq (D.E)^2/D^2$. 
\end{prop}
\begin{pf}
For a reference, see \cite[Proposition 2]{sakai1990reider}. 
\end{pf}

\begin{cor}\label{blowup}
Let $\pi: Y \to X$ be a birational morphism between two smooth surfaces and let $\wt{D}$ be a big Cartier divisor on $Y$ such that $\wt{D}^2>0$. 
Assume that $h^1(X,\cO_X(-\wt{D}))\neq 0$ and $\wt{D}$ is $m$-unstable for some $m>0$. 
That is, there exists a nonzero effective divisor $\wt{E}$ such that 
$m\wt{D} - 2\wt{E}$ is big and 
$(m\wt{D}-\wt{E}).\wt{E}\leq 0$.

Set $D = \pi_*\wt{D}$, $E = \pi_*\wt{E}$ and $\alpha = D^2-\wt{D}^2$. 
If $D$ is nef and $E$ is a nonzero effective divisor, then
$0\leq D.E<m\alpha/2$ and 
$mD.E-m^2\alpha/4\leq E^2\leq (D.E)^2/D^2$. 
\end{cor}
\begin{pf}
Write $\wt{B} = m\wt{D}$. 
Since $\wt{D}$ is $m$-unstable, $\wt{B}$ is 1-unstable. 
Thus, we can use Proposition~\ref{key_prop} above. 
Note that $\alpha_B = B^2-\wt{B}^2 = m^2(D^2-\wt{D}^2) = m^2\alpha_D$. \qed
\end{pf}

\section{Fujita's Conjecture for Quasi-Elliptic Surfaces}
From now on, $X$ and $Y$ are quasi-elliptic surfaces, and $A$ is an ample divisor on $X$. 
We first improve \cite[Proposition 4.3]{di2015effective}.

\begin{prop}\label{unstable}
Let $X$ be a quasi-elliptic surface with a quasi-elliptic fibration $f:X\rw C$ and $D$ be a big divisor on $X$ with $D^2>0$ and $h^1(X,\cO_X(-D))\neq 0$. 
Then $D$ is $p$-unstable. 

Moreover, let $F$ be a general fiber of the fibration $f$ and $E$ be a non-zero effective divisor whose existence is guaranteed by $p$-unstability of $D$ (see Definition~\ref{m_unstable}).
Then we have $(3D-2E).F = 1$ when $p=3$ and $(D-E).F = 1$ when $p=2$.
\end{prop}
\begin{pf}
Since $h^1(X,\cO_X(-D))\neq 0$, there exists a non-split short exact sequence 
\[\xymatrix{0\ar[r] & \cO_X \ar[r] & \cE \ar[r] & \cO_X(D) \ar[r] & 0}\]
given by a nonzero element of $\mbox{Ext}^1(\cO_X(D),\cO_X) \cong H^1(X,\cO_X(-D))$, 
where $\cE$ is a vector bundle of rank two. 
Note that $c_1^2(\cE)-4c_2(\cE) = D^2>0$. 
By Theorem~\ref{lift}, we have $(F^e)^*\cE$ is unstable for some $e>0$ and $\rho : Y \to X$ is a purely inseparable morphism of degree $p^e$.
Let \[\xymatrix{0\ar[r] & \cO_X(D_1)\ar[r] & F^{e*}\cE \ar[r] & \cI_Z\cdot\cO_X(D_2) \ar[r] & 0}\]
be a de-stabilizing sequence for $(F^e)^*\cE$. 

By Lemma~\ref{key_lem} and its proof, $D$ is $p^e$-unstable. 
Let $E$ be a non-zero effective divisor whose existence is guaranteed by $p^e$-unstability of $D$, $F$ be a general fiber of $f : X \to B$, and $C = \rho^*F$. 
Note that 
\begin{eqnarray*}
-K_Y.C &=& \rho^*\left(\frac{p^e-1}{p^e}(D_1-D_2)-K_X\right).C \\
&=& \rho^*\left(\frac{p^e-1}{p^e}(p^eD-2E)-K_X\right).C \\
&=& p^e\left(\frac{p^e-1}{p^e}(p^eD-2E)-K_X\right).F \\
&=& (p^e-1)(p^eD-2E).F >0
\end{eqnarray*}
where 
\begin{enumerate}
\item the first equality follows from Theorem~\ref{lift}, 
\item the second equality follows from 
\[D_1-D_2 = (D_1+D_2)-2D_2\equiv c_1(F^{e*}\cE)-2E\equiv p^eD-2E,\]
\item the third equality follows from projective formula, 
\item the fourth equality follows since $F$ has arithmetic genus one, and 
\item the last inequality follows since $p^eD-2E$ is big and $F$ is a general fiber of $f$. 
\end{enumerate}

Notice that $Y$ is a local complete intersection since $Y$ is a divisor in a smooth variety. 
Then by \cite[Corollary 2.14]{di2015effective}, we have $-K_Y.C\leq 3$. 
This gives 
\[3\geq -K_Y.C = (p^e-1)(p^eD-2E).F.\]

When $p=3$, we have $(p^e-1)(p^eD-2E).F\geq 3^e-1\geq 8$ if $e\geq 2$, which is impossible. 

When $p=2$, we have $(p^e-1)(p^eD-2E).F = 2(2^e-1)(2^{e-1}D-E).F\geq 2(2^e-1)\geq 6$ if $e\geq 2$, which is impossible. 

Thus, $e$ must be $1$ and $D$ is $p$-unstable.  

Moreover, when $p=3$, we have 
\[(p^e-1)(p^eD-2E).F = 2(3D-2E).F\]
which is a positive even integer less than $3$. 
So we have $(3D-2E).F=1$. 
When $p=2$, we have 
\[(p^e-1)(p^eD-2E).F = 2(D-E).F\] 
which is a positive even integer less than $3$. 
So we have $(D-E).F=1$. 
\qed
\end{pf}

Now we are ready to prove 
\begin{thm}\label{main}
Fujita's conjecture is true for quasi-elliptic surfaces $X$. 
That is, given a quasi-elliptic surface $X$ and any ample line bundle $A$ on $X$, we have 
\begin{enumerate}
\item $K_X+tA$ is base point free for $t\geq 3$; and 
\item $K_X+tA$ is very ample for $t\geq 4$. 
\end{enumerate}
\end{thm}

\begin{pf}
We divide the proof into several steps. 
We first prove that $K_X+tA$ is base point free for $t\geq 3$.

\begin{enumerate}
\item[(Step 1)] (Preparation.) Let $D = tA$ and assume that $|K_X+D|$ has a base point at $x\in X$. 
Let $\pi : Y \to X$ be the blow-up at $x$. 
Since $x$ is a base point, we have that 
\[h^1(Y,\cO_Y(K_Y+\pi^*D-2E_x)) = h^1(X,\cO_X(K_X+D)\otimes \mathfrak{m}_x) \neq 0\] where $E_x$ is the exceptional divisor of $\pi$. 
Let $\wt{D} = \pi^*D-2E_x$. 

In order to apply Proposition~\ref{unstable}, we need to check that $\wt{D}$ is big and $\wt{D}^2>0$. 
Note that 
\begin{eqnarray*}
h^0(Y,\cO_Y(\ell\wt{D})) &=& h^0(Y,\cO_Y(\ell(\pi^*D-2E_x))) \\
&=& h^0(X, \cO_X(\ell D)\otimes \mathfrak{m}_x^{2\ell}) \\
&\geq & \frac{D^2}{2}\ell^2 + O(\ell) - 
{\left(\begin{matrix}
2\ell + 1\\
2
\end{matrix}\right)} \\
&=& \frac{t^2A^2-4}{2}\ell^2+O(\ell).
\end{eqnarray*}
So $\wt{D}$ is big whenever $t\geq 3$. 
Also note that $\wt{D}^2 = D^2-4 = t^2A^2-4\geq 5$ when $t\geq 3$. 

Applying Proposition~\ref{unstable} on $Y$ and $\wt{D}$, we have that $\wt{D}$ is $p$-unstable. 
So there is a nonzero effective divisor $\wt{E}$ such that $p\wt{D}-2\wt{E}$ is big and $(p\wt{D}-\wt{E}).\wt{E}\leq 0$. 
Let $E = \pi_*\wt{E}$. Note that $E$ is an effective divisor. 
If $E=0$, then $\wt{E} = bE_x$ for $b>0$. 
Thus, $(p\wt{D}-\wt{E}).\wt{E} = b(b+2p)>0$, which is a contradiction. 
Therefore, $E$ is non-zero. 

Also $\pi_*\wt{D} = D = tA$ is ample and $\alpha = D^2 - \wt{D}^2 = 4$. 
Hence, by Corollary~\ref{blowup}, we have 
\begin{equation}\label{bpf_ineq}
0< tA.E < 2p\leq 6 \mbox{ and } 
ptA.E-p^2\leq E^2\leq (A.E)^2/A^2.
\end{equation}
So we have $0<A.E < \frac{6}{t}\leq 2$ and thus, $A.E = 1$ and $E$ is an irreducible curve. 
The second inequality in~(\ref{bpf_ineq}) becomes 
\begin{equation}\label{bpf}
pt-p^2\leq E^2\leq 1/A^2\leq 1.
\end{equation}

\item[(Step 2)] If $p=2$, then $2\leq 2t-4 \leq E^2\leq 1$, which is impossible. 

\item[(Step 3)] If $p=3$, then $3t-9\leq E^2\leq 1$. 
This happens only when $t=3$ and $E^2=0$ or $1$. 

Now, by Proposition~\ref{unstable}, we have 
\begin{equation}\label{p3bpf}
1=(3\wt{D}-2\wt{E}).\pi^*F = (9A-2E).F. 
\end{equation}
Since $F$ is nef and $A$ is ample, we get $9A.F\geq 9$ and so, by equality~(\ref{p3bpf}), we have 
\begin{equation}\label{cf4}
E.F\geq 4. 
\end{equation} 
Note that $E+F$ is nef since $(E+F).E \geq 0 + F.E\geq 4$ and $(E+F).F=E.F\geq 4$. 

\item[(Step 4)] If $E^2 = 1$, then $A^2 = 1$ by inequality~(\ref{bpf}) and $A$ is numerically equivalent to $E$ by Hodge inequality. 
Thus, by equality~(\ref{p3bpf}), we have $7A.F = 1$, which is impossible. 

\item[(Step 5)] So we have $E^2 = 0$. 
Applying Lemma~\ref{Hodeg_ineq} to $9A-2E$ and $E+F$, we have 
\[(9A-2E)^2(E+F)^2\leq \left((9A-2E).(E+F)\right)^2.\]
Thus, we have 
\[(81A^2-36)(2F.E)\leq (9A.E+(9A-2E).F)^2 = 100\] 
since $A.E=1$ and $(9A-2E).F = 1$ by equality~(\ref{p3bpf}). 
Thus, 
\[5\leq 9A^2-4\leq \frac{100}{18(F.E)} \leq \frac{100}{18\times 4}\leq 2, \mbox{ which is impossible.}\]

Hence, we have shown that $K_X+tA$ is base point free whenever $t\geq 3$. 
\end{enumerate}

Now we prove $K_X+tA$ is very ample when $t\geq 4$. 
\begin{enumerate}
\item[(Step 1)] (Preparation.) Let $D = tA$. 
It suffices to show $|K_X+D|$ separates points and tangents. 
(For a reference, see~\cite[Proposition II.7.3]{hartshorne1977algebraic}.) 
Assume that $|K_X+D|$ does not separate points $x$ and $y$ (resp. does not separate tangents at $x$). 
Then we have 
\[\begin{array}{l}
h^1(Y,\cO_Y(K_Y+\pi^*D-2E_x-2E_y)) = h^1(X,\cO_X(K_X+D)\otimes\mathfrak{m}_x\otimes\mathfrak{m}_y) \neq 0 \\
(\mbox{resp. } h^1(Y,\cO_Y(K_Y+\pi^*D-3E_x)) = h^1(X,\cO_X(K_X+D)\otimes\mathfrak{m}^2_x) \neq 0)
\end{array}\]
where $\pi : Y\to X$ is the blow-up of $X$ at $x, y$ and $E_x, E_y$ are the exceptional divisor 
(resp. $\pi : Y\to X$ is the blow-up of $X$ at $x$ and $E_x$ is the exceptional divisor.)

Now let $\wt{D} = \pi^*D-2E_x-2E_y$ (resp. $\wt{D} = \pi^*D-3E_x$).
By the similar argument as above, $\wt{D}$ is big and $\wt{D}^2>0$ whenever $t\geq 4$. 

Applying Proposition~\ref{unstable} to $Y$ and $\wt{D}$, we have that $\wt{D}$ is $p$-unstable. 
So there is a nonzero effective divisor $\wt{E}$ such that $p\wt{D}-2\wt{E}$ is big and $(p\wt{D}-\wt{E}).\wt{E}\leq 0$. 
Let $E = \pi_*\wt{E}$. Note that $E$ is a non-zero effective divisor by the similar argument as above. 
Also we have $\pi_*\wt{D} = D = tA$ is ample and $\alpha = D^2 - \wt{D}^2 = 8$ (resp. $9$). 
Hence, by Corollary~\ref{blowup}, we have 
\begin{equation}\label{va}
0< tA.E < p\alpha/2 \mbox{ and }
ptA.E-p^2\alpha/4\leq E^2\leq (A.E)^2/A^2.
\end{equation}

\item[(Step 2)] When $p=3$, by Proposition~\ref{unstable}, we have 
\[(3\wt{D}-2\wt{E}).\pi^*F = 1\]
and thus
\begin{equation}\label{va2}
(3tA-2E).F=1.
\end{equation}
If $t$ is even, then the left hand side is an even integer, which is impossible. 

\item[(Step 3)] So $t$ is odd. 
From inequalities~(\ref{va}), we have  
\begin{equation}\label{vap3kodd}
0< tA.E < \frac{3}{2}\alpha\leq \frac{27}{2} \mbox{ and }
3tA.E-\frac{9}{4}\alpha\leq E^2\leq (A.E)^2/A^2.
\end{equation}

Then we have $A.E=1$ or $2$ since $A.E<\frac{27}{2t}\leq\frac{27}{10}<3$. 
If $A.E=2$, then we have 
\[9< 6t-\frac{81}{4}\leq 6t - \frac{9}{4}\alpha \leq E^2\leq 4/A^2\leq 4, \mbox{ which is impossible.}\]
Thus $A.E=1$ and so $E$ is an irreducible curve. 
Also, from the second inequality in~(\ref{vap3kodd}), we have 
\[3t-\frac{81}{4}\leq 3t-\frac{9}{4}\alpha\leq E^2\leq\frac{1}{A^2}\leq 1.\] 
Hence $t$ is 5 or 7 and 
\begin{equation}\label{E2geq-5}
-5\leq E^2\leq 1.
\end{equation}
Using equality~(\ref{va2}), we have $1+2E.F = 3tA.F\geq 15$. 
So 
\begin{equation}\label{efgeq7}
E.F\geq 7
\end{equation} 
and $E+F$ is nef since $(E+F).E\geq -5+7=2$ and $(E+F).F\geq 7$. 

Applying Lemma~\ref{Hodeg_ineq} to $3tA-2E$ and $E+F$, we have  
\[(3tA-2E)^2(E+F)^2\leq ((3tA-2E).(E+F))^2.\]
Note that the left hand side
\begin{eqnarray*}
(3tA-2E)^2(E+F)^2 &=& (9t^2A^2-12t+4E^2)(E^2+2E.F) \\
&\geq & (4E^2+141)(E^2+2E.F) \\
&\geq & (4E^2+141)(E^2+14)
\end{eqnarray*}
where 
\begin{enumerate}
\item the first inequality follows from $A^2\geq 1$, $5\leq t\leq 7$, and nefness of $E+F$; and
\item the second inequality follows from inequalities~(\ref{E2geq-5}) and (\ref{efgeq7}). 
\end{enumerate}
Also the right hand side
\begin{eqnarray*}
((3tA-2E).(E+F))^2 &=& (1+3t-2E^2)^2 \\
&\leq & (22-2E^2)^2 
\end{eqnarray*}
where 
\begin{enumerate}
\item the equality follows from equality~(\ref{va2}); and 
\item the inequality follows from inequality (\ref{E2geq-5}) and $t\leq 7$. 
\end{enumerate}
To sum up, we have $(4E^2+141)(E^2+14)\leq (22-2E^2)^2$ 
and thus $E^2\leq -6$, which contradicts to inequality~(\ref{E2geq-5}). 

\item[(Step 4)] Now we deal with $p=2$. 
The inequalities~(\ref{va}) becomes 
\begin{equation}\label{vap2}
0< tA.E < \alpha\leq 9 \mbox{ and }
2tA.E-\alpha\leq E^2\leq (A.E)^2/A^2.
\end{equation}
Hence, $A.E = 1$ or $2$ since $A.E<\frac{9}{t}\leq \frac{9}{4}$. 
If $A.E = 2$, then we have $7\leq 4t-\alpha \leq E^2\leq 4/A^2\leq 4$, which is impossible. 

\item[(Step 5)] Thus, we have $A.E = 1$, and so $E$ is an irreducible curve. 
Then, from the second inequality in~(\ref{vap2}), we have 
\[2t-9\leq 2t-\alpha\leq E^2\leq 1/A^2\leq 1.\] 
So $t=5$ and $E^2=1$; or $t=4$ and $-1\leq E^2\leq 1$.
By Proposition~\ref{unstable}, we have $(\wt{D}-\wt{E}).\pi^*F=1$ and thus 
\begin{equation}\label{va10}
(tA-E).F=1.
\end{equation}
So we have 
\begin{equation}\label{va11}
E.F = tA.F-1\geq 3. 
\end{equation}
Therefore, $E+F$ is nef since $(E+F).E\geq -1+3=2$ and $(E+F).F\geq 3$. 

\item[(Step 6)] If $t = 5$, then we have $E^2 = 1$ and thus $A\equiv E$ by Lemma~\ref{Hodeg_ineq}. 
But, from equality (\ref{va10}), we have \[1 = (5A-E).F = 4A.F, \mbox{ which is impossible.}\] 

\item[(Step 7)] Thus $t=4$. 
Applying Lemma~\ref{Hodeg_ineq} to $4A-E$ and $E+F$, we have 
\begin{equation}\label{va_p2}
(4A-E)^2(E+F)^2\leq \left((4A-E).(E+F)\right)^2 = (1+4-E^2)^2=(5-E^2)^2.
\end{equation}

\item[(Step 8)] When $E^2 = 1$, we have $A\equiv E$ by Lemma~\ref{Hodeg_ineq}. 
Thus, from equality (\ref{va10}), we have \[1=(4A-E).F = 3A.F, \mbox{ which is impossible.}\] 

\item[(Step 9)] When $E^2 = 0$, we have, from inequality~(\ref{va_p2}),
\[(16A^2-8)(2E.F)\leq 25\]
and thus, by inequality (\ref{va11}),  
\[1\leq 2A^2-1\leq \frac{25}{16(E.F)}\leq \frac{25}{16\times 3}<1, \mbox{ which is impossible.}\]

\item[(Step 10)] When $E^2 = -1$, from inequality~(\ref{va_p2}), we have 
\[(16A^2-9)(2E.F-1)\leq 36\]
and thus, by inequality (\ref{va11}), 
\[7\leq 16A^2-9\leq \frac{36}{2E.F-1}\leq \frac{36}{5}<8.\]
Therefore, $A^2 = 1$ and $E.F=3$. 
Moreover, by equality (\ref{va10}), we have $A.F=1$. 
Now applying Lemma~\ref{Hodeg_ineq} to $A$ and $E+F$, we have 
\[5 = A^2(E+F)^2\leq (A.(E+F))^2 = 4, \mbox{ which is impossible.}\]

Hence, we have shown that $K_X+tA$ is very ample whenever $t\geq 4$. 
\end{enumerate}
\qed
\end{pf}

\bibliographystyle{amsalpha}
\addcontentsline{toc}{chapter}{\bibname}
\normalem
\bibliography{FujitaConjectureQuasiEllipticSurfaces}

\end{document}